\documentclass[reqno]{amsart}
\usepackage{amssymb}
\usepackage{ifpdf}
\ifpdf
 \usepackage[hyperindex,pagebackref]{hyperref}%
\else
 \expandafter\ifx\csname dvipdfm\endcsname\relax
 \usepackage[hypertex,hyperindex,pagebackref]{hyperref}
 \else
 \usepackage[dvipdfm,hyperindex,pagebackref]{hyperref}
 \fi
\fi
\allowdisplaybreaks[4]
\numberwithin{equation}{section}
\numberwithin{table}{section}
\theoremstyle{plain}
\newtheorem{thm}{Theorem}[section]

\begin{document}

\title[A new formula for Bernoulli and Genocchi numbers]
{A new explicit formula for Bernoulli and Genocchi numbers in terms of Stirling numbers}

\author[B.-N. Guo]{Bai-Ni Guo}
\address[Guo]{School of Mathematics and Informatics, Henan Polytechnic University, Jiaozuo City, Henan Province, 454010, China}
\email{\href{mailto: B.-N. Guo <bai.ni.guo@gmail.com>}{bai.ni.guo@gmail.com}, \href{mailto: B.-N. Guo <bai.ni.guo@hotmail.com>}{bai.ni.guo@hotmail.com}}
\urladdr{\url{http://www.researchgate.net/profile/Bai-Ni_Guo}}

\author[F. Qi]{Feng Qi}
\address[Qi]{College of Mathematics, Inner Mongolia University for Nationalities, Tongliao City, Inner Mongolia Autonomous Region, 028043, China; Department of Mathematics, College of Science, Tianjin Polytechnic University, Tianjin City, 300387, China}
\email{\href{mailto: F. Qi <qifeng618@gmail.com>}{qifeng618@gmail.com}, \href{mailto: F. Qi <qifeng618@hotmail.com>}{qifeng618@hotmail.com}, \href{mailto: F. Qi <qifeng618@qq.com>}{qifeng618@qq.com}}
\urladdr{\url{http://qifeng618.wordpress.com}}

\begin{abstract}
In the paper, the authors review some explicit formulas and establish a new explicit formula for Bernoulli and Genocchi numbers in terms of Stirling numbers of the second kind.
\end{abstract}



\thanks{This paper was typeset using \AmS-\LaTeX}

\maketitle

\section{Main results}

It is well known that Bernoulli numbers $B_{n}$ for $n\ge0$ may be defined by the power series expansion
\begin{equation}
\frac{x}{e^x-1}=\sum_{n=0}^\infty B_n\frac{x^n}{n!}=1-\frac{x}2+\sum_{k=1}^\infty B_{2k}\frac{x^{2k}}{(2k)!}, \quad \vert x\vert<2\pi,
\end{equation}
that Euler polynomials $E_n(x)$ are defined by
\begin{equation}
\frac{2e^{xt}}{e^t+1}=\sum_{n=0}^\infty E_n(x)\frac{t^n}{n!},
\end{equation}
that Genocchi numbers $G_n$ for $n\in\mathbb{N}$ are given by the generating function
\begin{equation}\label{Genocchi-dfn=eq}
\frac{2t}{e^t+1}=\sum_{n=1}^\infty G_n\frac{t^n}{n!},
\end{equation}
and that Stirling numbers of the second kind which may be generated by
\begin{equation}\label{2stirling-gen-funct-exp}
\frac{(e^x-1)^k}{k!}=\sum_{n=k}^\infty S(n,k)\frac{x^n}{n!}, \quad k\in\mathbb{N}
\end{equation}
and may be computed by
\begin{equation}\label{Stirling-Number-dfn}
S(k,m)=\frac1{m!}\sum_{\ell=1}^m(-1)^{m-\ell}\binom{m}{\ell}\ell^{k}, \quad 1\le m\le k.
\end{equation}
By the way, Stirling number of the second kind $S(n,k)$ may be interpreted combinatorially as the number of ways of partitioning a set of $n$ elements into $k$ nonempty subsets.
\par
Bernoulli numbers $B_n$ for $n\in\{0\}\cup\mathbb{N}$ satisfy
\begin{equation}
B_0=1,\quad B_1=-\frac12, \quad B_{2n+2}\ne0, \quad B_{2n+3}=0.
\end{equation}
For $n\in\mathbb{N}$, Genocchi numbers meet $G_{2n+1}=0$. The first few Genocchi numbers $G_{n}$ are listed in Table~\ref{Genocchi-Table}.
\begin{table}[htbp]
\caption{The first few Genocchi numbers $G_{n}$}\label{Genocchi-Table}
\begin{tabular}{|c|c|c|c|c|c|c|c|c|c|c|}
  \hline
  $n$   &$ 1 $&$ 2 $&$ 4 $&$ 6 $&$ 8 $&$ 10 $&$ 12 $&$ 14 $&$ 16 $&$ 18 $\\  \hline
  $G_n$ &$ 1 $&$ -1 $&$ 1 $&$ -3 $&$ 17 $&$ -155 $&$ 2073 $&$ -38227 $&$ 929569 $&$ -28820618 $\\
  \hline
\end{tabular}
\end{table}
Genocchi numbers $G_{2n}$ may be represented in terms of Bernoulli numbers $B_{2n}$ and Euler polynomials $E_{2n-1}(0)$ as
\begin{equation}
G_{2n}=2(1-2^{2n})B_{2n}=2nE_{2n-1}(0), \quad n\in\mathbb{N}.
\end{equation}
See~\cite[p.~49]{Comtet-Combinatorics-74}. As a result, we have
\begin{equation}
G_{n}=2(1-2^{n})B_{n}, \quad n\in\mathbb{N}.
\end{equation}
\par
The first formula for Bernoulli numbers $B_n$ listed in~\cite{Gould-MAA-1972} is
\begin{equation}\label{Higgins-Gould-B}
B_n=\sum_{k=0}^n\frac1{k+1}\sum_{j=0}^k(-1)^j\binom{k}{j}j^n,\quad n\ge0,
\end{equation}
which is a special case of the general formula~\cite[(2.5)]{Higgins-JLMS-1970}. The formula~\eqref{Higgins-Gould-B} is equivalent to
\begin{equation}\label{Bernoulli-Stirling-eq}
B_n=\sum_{k=0}^n(-1)^k\frac{k!}{k+1}S(n,k), \quad n\in\{0\}\cup\mathbb{N},
\end{equation}
which was listed in~\cite[p.~536]{GKP-Concrete-Math-1st} and~\cite[p.~560]{GKP-Concrete-Math-2nd}. Recently, four alternative proofs of the formula~\eqref{Bernoulli-Stirling-eq} were provided in~\cite{ANLY-D-12-1238.tex, Bernoulli-Stirling-4P.tex}. A generalization of the formula~\eqref{Bernoulli-Stirling-eq} was supplied in~\cite{BernPolRep.tex}. In all, we may collect at least seven alternative proofs for the formula~\eqref{Higgins-Gould-B} or~\eqref{Bernoulli-Stirling-eq} in~\cite{Gould-MAA-1972, GKP-Concrete-Math-2nd, ANLY-D-12-1238.tex, Higgins-JLMS-1970, Logan-Bell-87, Bernoulli-Stirling-4P.tex} and closely related references therein.
\par
In~\cite[p.~48, (11)]{Gould-MAA-1972}, it was deduced that
\begin{equation}\label{Higgins-Gould-B(11)}
B_n=\sum_{j=0}^n(-1)^j\binom{n+1}{j+1}\frac{n!}{(n+j)!}\sum_{k=0}^j(-1)^{j-k}\binom{j}{k}k^{n+j}, \quad n\ge0,
\end{equation}
which may be rearranged as
\begin{equation}\label{Bernoulli-Stirling-formula}
B_n=\sum_{i=0}^n(-1)^{i}\frac{\binom{n+1}{i+1}}{\binom{n+i}{i}}S(n+i,i), \quad n\ge0.
\end{equation}
The formula~\eqref{Bernoulli-Stirling-formula} was rediscovered in the preprint~\cite{Bernoulli-No-Int-New.tex}.
On 21 January 2014, the authors searched out that the formula~\eqref{Bernoulli-Stirling-formula} was also derived in~\cite[p.~59]{Jeong-Kim-Son-JNT-2005} and~\cite[p.~140]{Shirai-Sato-JNT-2001}.
\par
In~\cite[p.~1128, Corollary]{recursion}, among other things, it was found that
\begin{equation}\label{Bernoulli-N-Guo-Qi-99}
B_{2k}= \frac12 - \frac1{2k+1} - 2k \sum_{i=1}^{k-1}
\frac{A_{2(k-i)}}{2(k - i) + 1}
\end{equation}
for $k\in\mathbb{N}$, where $A_m$ is defined by
\begin{equation*}
\sum_{m=1}^nm^k=\sum_{m=0}^{k+1}A_mn^{m}.
\end{equation*}
\par
In~\cite[Theorem~1.4]{Tan-Cot-Bernulli-No.tex}, among other things, it was presented that
\begin{equation}\label{Bernoulli-N-Explicit}
B_{2k}=\frac{(-1)^{k-1}k}{2^{2(k-1)}(2^{2k}-1)}\sum_{i=0}^{k-1}\sum_{\ell=0}^{k-i-1} (-1)^{i+\ell}\binom{2k}{\ell}(k-i-\ell)^{2k-1}, \quad k\in\mathbb{N}.
\end{equation}
In~\cite[Theorem~3.1]{exp-derivative-sum-Combined.tex}, it was obtained that
\begin{multline}\label{Bernumber-formula-eq}
B_{2k}=1+\sum_{m=1}^{2k-1}\frac{S(2k+1,m+1) S(2k,2k-m)}{\binom{2k}{m}} \\*
-\frac{2k}{2k+1}\sum_{m=1}^{2k}\frac{S(2k,m)S(2k+1,2k-m+1)}{\binom{2k}{m-1}}, \quad k\in\mathbb{N}.
\end{multline}
\par
The aim of this paper is to find the following new explicit formula for Bernoulli numbers $B_k$, or say, Genocchi numbers $G_{k}$, in terms of Stirling numbers of the second kind $S(k,m)$.

\begin{thm}\label{Genocchi-Stirling=thm}
For all $k\in\mathbb{N}$, Genocchi numbers $G_{k}$ may be computed by
\begin{equation}
G_k=2(1-2^{k})B_{k}=(-1)^{k}k\sum_{m=1}^{k}(-1)^m\frac{(m-1)!}{2^{m-1}}S(k,m).
\end{equation}
\end{thm}

\section{Proof of Theorem~\ref{Genocchi-Stirling=thm}}

Differentiating on both sides of the equation~\eqref{Genocchi-dfn=eq} and employing Leibniz identity for differentiation give
\begin{gather*}
\biggl(\frac{2t}{e^t+1}\biggr)^{(k)}
=2\biggl[t\biggl(\frac1{e^t+1}\biggr)^{(k)}
+k\biggl(\frac1{e^t+1}\biggr)^{(k-1)}\biggr]
=\sum_{n=k}^\infty G_n\frac{t^{n-k}}{(n-k)!}.
\end{gather*}
In~\cite[Theorem~2.1]{Eight-Identy-More.tex} and~\cite[Theorem~3.1]{CAM-D-13-01430-Xu-Cen}, it was obtained that, when $\lambda>0$ and $t\ne-\frac{\ln\lambda}\alpha$ or when $\lambda<0$ and $t\in\mathbb{R}$,
\begin{equation}\label{id-gen-new-form1}
\biggl(\frac1{\lambda e^{\alpha t}-1}\biggr)^{(k)}
=(-1)^k\alpha^k\sum_{m=1}^{k+1}{(m-1)!S(k+1,m)}\biggl(\frac1{\lambda e^{\alpha t}-1}\biggr)^m.
\end{equation}
Specially, when $\lambda=-1$ and $\alpha=1$, the identity~\eqref{id-gen-new-form1} becomes
\begin{equation}\label{alpha=1lambda=-1}
\biggl(\frac1{e^{t}+1}\biggr)^{(k)}
=(-1)^{k+1}\sum_{m=1}^{k+1}(-1)^m{(m-1)!S(k+1,m)}\biggl(\frac1{e^{t}+1}\biggr)^m.
\end{equation}
Consequently, it follows that
\begin{gather*}
G_k=\lim_{t\to0}\sum_{n=k}^\infty G_n\frac{t^{n-k}}{(n-k)!}
=2k\lim_{t\to0}\biggl(\frac1{e^t+1}\biggr)^{(k-1)}\\
=2k(-1)^{k}\sum_{m=1}^{k}(-1)^m{(m-1)!S(k,m)}\lim_{t\to0}\biggl(\frac1{e^{t}+1}\biggr)^m\\
=(-1)^{k}k\sum_{m=1}^{k}(-1)^m\frac{(m-1)!}{2^{m-1}}S(k,m).
\end{gather*}
The proof of Theorem~\ref{Genocchi-Stirling=thm} is complete.

\end{document}